\documentclass{amsart}

\usepackage{graphicx}

\numberwithin{equation}{section}

\def\be{\begin{equation}}
\def\ee{\end{equation}}

\newcommand{\bthm}{\begin{thm}}
\newcommand{\ethm}{\end{thm}}

\newcommand{\beqq}{\begin{eqnarray*}}
\newcommand{\eeq}{\end{eqnarray}}
\newcommand{\eeqq}{\end{eqnarray*}}

\newcommand{\A}{\mathcal A}

\newcommand{\IC}{{\mathbb C}}
\newcommand{\ID}{{\mathbb D}}

\newtheorem{thm}{Theorem}

\theoremstyle{definition}
\newtheorem{defn}{Definition}

\title[\uppercase{Hankel determinant for a class of analytic  functions}]{\uppercase{Hankel determinant for a class of analytic  functions}}

\author[M. Obradovi\'{c}]{Milutin Obradovi\'{c}}
\address{Department of Mathematics\newline \indent
Faculty of Civil Engineering\newline \indent
University of Belgrade\newline \indent
Bulevar Kralja Aleksandra 73\newline \indent
11000, Belgrade, Serbia}
\email{obrad@grf.bg.ac.rs}

\author[N. Tuneski]{Nikola Tuneski}
\address{Department of Mathematics and Informatics
\newline \indent
Faculty of Mechanical Engineering
\newline \indent
Ss. Cyril and Methodius University in Skopje
\newline \indent
Karpo\v{s} II b.b., 1000 Skopje
\newline \indent
Republic of North Macedonia}
\email{nikola.tuneski@mf.edu.mk}

\subjclass[2000]{30C45, 30C50}

\keywords{analytic, univalent, Hankel determinant.}

\begin{document}

\begin{abstract}
Let $f$ be analutic in the unit disk $\ID$ and normalized so that $f(z)=z+a_2z^2+a_3z^3+\cdots$. In this paper we give sharp bound of Hankel determinant of the second order  for the class of analytic unctions satisfying
\[ \left|\arg \left[\left(\frac{z}{f(z)}\right)^{1+\alpha}f'(z) \right] \right|<\gamma\frac{\pi}{2}   \quad\quad (z\in\ID),\]
for $0<\alpha<1$ and $0<\gamma\leq1$.
\end{abstract}

\maketitle

\section{Introduction and preliminaries}


Let ${\mathcal A}$ denote the family of all analytic functions
in the unit disk $\ID := \{ z\in \IC:\, |z| < 1 \}$ and  satisfying the normalization
$f(0)=0= f'(0)-1$.

A function $f\in \mathcal{A}$ is said to be {\it strongly starlike of order $\beta, 0<\beta \leq 1$} if, and only if,
$$\left|\arg\frac{zf'(z)}{f(z)}\right|<\beta\frac{\pi}{2}\quad\quad (z\in\ID).$$
We denote this class by $\mathcal{S}^{\star}_{\beta}$. If $\beta=1$,
then $\mathcal{S}^{\star}_{1}\equiv \mathcal{S}^{\star}$ is the well-known class of \textit{starlike functions}.

In \cite{MO-1998} the author introduced the class $\mathcal{U}(\alpha, \lambda)$ ($0<\alpha$ and $\lambda <1$)
consisting of functions $f\in \mathcal{A}$  for which we have
$$
\left|\left(\frac{z}{f(z)}\right)^{1+\alpha}f'(z)-1\right|<\lambda \quad\quad (z\in\ID).
$$
In the same paper it is shown that $\mathcal{U}(\alpha, \lambda)\subset \mathcal{S}^{\star}$ if
$$0<\lambda \leq \frac{1-\alpha}{\sqrt{(1-\alpha)^{2}+\alpha^{2}}}.$$
The most valuable up to date results about this class can be found in Chapter 12 from \cite{DTV-book}.

In the paper \cite{MO-2000}  the author considered univalence of the class of functions $f\in \mathcal{A}$
satisfying the condition
\be\label{eq-17}
 \left|\arg \left[\left(\frac{z}{f(z)}\right)^{1+\alpha}f'(z) \right] \right|<\gamma\frac{\pi}{2}   \quad\quad (z\in\ID)
\ee
for $0<\alpha<1$ and $0<\gamma\leq1$, and proved the following theorem.

\medskip

\noindent {\it {\bf Theorem A}. Let $f\in \mathcal{A}$, $0<\alpha<\frac{2}{\pi}$ and let
$$ \left|\arg \left[\left(\frac{z}{f(z)}\right)^{1+\alpha}f'(z) \right]\right|<\gamma_{\star}(\alpha)\frac{\pi}{2} \quad\quad (z\in\ID),$$
where
$$ \gamma_{\star}(\alpha)=\frac{2}{\pi}\arctan \left(\sqrt{\frac{2}{\pi\alpha}-1}\right)-\alpha \sqrt{\frac{2}{\pi\alpha}-1}.$$
Then $f\in\mathcal{S}^{\star}_{\beta},$ where
$$\beta=\frac{2}{\pi}\arctan \sqrt{\frac{2}{\pi\alpha}-1}.$$ }

\section{Main result}

In this paper we will give the sharp estimate for Hankel determinant of the second order  for the class of analytic unctions $f\in \mathcal{A}$
which satisfied the condition \eqref{eq-17}.

\begin{defn}\label{def-hank}
Let $f\in \A$. Then the $qth$ Hankel determinant of $f$ is defined for $q\geq 1$, and
$n\geq 1$ by
\[
        H_{q}(n) = \left |
        \begin{array}{cccc}\nonumber 
        a_{n} & a_{n+1}& \ldots& a_{n+q-1}\\
        a_{n+1}&a_{n+2}& \ldots& a_{n+q}\\
        \vdots&\vdots&~&\vdots \\
        a_{n+q-1}& a_{n+q}&\ldots&a_{n+2q-2}\\
        \end{array}
        \right |.
\]
\end{defn}
Thus, the second Hankel determinant is $H_{2}(2)= a_2a_4-a_{3}^2$.

\bthm\label{18-th-2}
Let $f(z)=z+a_{2}z^{2}+a_{3}z^{3}+\cdots$ belongs to the class $\mathcal{A}$ and satisfy the condition \eqref{eq-17}. Then we have the  next sharp estimation:
$$  |H_{2}(2)|=|a_{2}a_{4}-a_{3}^{2}|\leq \left(\frac{2\gamma}{2-\alpha}\right)^{2},$$
where $0<\alpha< 2-\sqrt{2}$ and $0<\gamma \leq\frac{1}{2}(\alpha^{2}-4\alpha +2).$
\ethm

\begin{proof}
We can write the condition \eqref{eq-17}  in the form
\be\label{eq-18-1}
\left(\frac{f(z)}{z}\right)^{-(1+\alpha)}f'(z)=\left(\frac{1+\omega(z)}{1-\omega (z)}\right)^{\gamma}
\,\left (=(1+2\omega(z)+2\omega^{2}(z)+\cdots)^{\gamma}\right),
\ee
where $\omega$ is analytic in $\ID$ with $\omega(0)=0$ and $|\omega(z)|<1$, $z\in\ID$.
If we denote by $L$ and $R$ left and right hand side of equality \eqref{eq-18-1}, then we have
\beqq
\begin{split}
L&= \left[1-(1+\alpha)(a_{2}z+\cdots)+{-(1+\alpha)\choose 2}(a_{2}z+\cdots)^{2} \right.\\
&+\left.{-(1+\alpha)\choose3}(a_{2}z+\cdots)^{3}+\cdots\right] \cdot(1+2a_{2}z+3a_{3}z^{2}+4a_{4}z^{3}+\cdots)
\end{split}
\eeqq
and if we put $\omega(z)=c_{1}z+c_{2}z^{2}+\cdots$ :
\beqq
\begin{split}
R&= 1+\gamma\left[2(c_{1}z+c_{2}z^{2}+\cdots)+2(c_{1}z+c_{2}z^{2}+\cdots)^{2}+\cdots\right]\\
&+ {\gamma \choose 2}\left[2(c_{1}z+c_{2}z^{2}+\cdots)+2(c_{1}z+c_{2}z^{2}+\cdots)^{2}+\cdots\right]^{2}\\
&+ {\gamma \choose3} \left[2(c_{1}z+c_{2}z^{2}+\cdots)+2(c_{1}z+c_{2}z^{2}+\cdots)^{2}+\cdots\right]^{3}+\cdots.
\end{split}
\eeqq
If we compare the coefficients on $z,z^{2},z^{3}$ in $L$ and $R$, then, after some calculations, we obtain
\be\label{eq-19}
\begin{array}{l}
\displaystyle\smallskip
a_{2}=\frac{2\gamma}{1-\alpha}c_{1}, \\
\displaystyle\smallskip
a_{3}=\frac{2\gamma}{2-\alpha}c_{2}+\frac{2(3-\alpha)\gamma^{2}}{(1-\alpha)^{2}(2-\alpha)}c_{1}^{2},\\
\displaystyle\smallskip
a_{4}=\frac{2\gamma}{3-\alpha}\left(c_{3}+\mu c_{1}c_{2}+\nu c_{1}^{3}\right),
\end{array}
\ee
where
\be\label{eq-20}
\mu=\mu (\alpha,\gamma)=\frac{2(5-\alpha)\gamma}{(1-\alpha)(2-\alpha)}\quad\mbox{and}\quad
\nu=\nu (\alpha, \gamma)=\frac{1}{3}+
\frac{2}{3}\frac{(\alpha ^{2}-6\alpha +17)\gamma^{2}}{(1-\alpha)^{3}(2-\alpha )}.
\ee
By using the relations \eqref{eq-19} and \eqref{eq-20}, after some simple computations, we obtain
\[
H_{2}(2)=\frac{4\gamma^{2}}{(1-\alpha)(3-\alpha)}\left(c_{1}c_{3}+\mu_{1}c_{1}^{2}c_{2}+(\frac{1}{3}-\nu_{1})
c_{1}^{4}-\frac{(1-\alpha)(3-\alpha)}{(2-\alpha)^{2}}c_{2}^{2}\right),
\]
where
\[
\mu_{1}=\frac{2\gamma}{(2-\alpha)^{2}},\quad
\nu_{1}=\frac{(\alpha^{2}-10\alpha +13)\gamma^{2}}{3(1-\alpha)^{2}(2-\alpha)^{2}},
\]
and from here
\be\label{eq-23}
\begin{split}
|H_{2}(2)| &\leq \frac{4\gamma^{2}}{(1-\alpha)(3-\alpha)}\bigg(|c_{1}||c_{3}|+\mu_{1}|c_{1}|^{2}|c_{2}|\\
&+\left|\frac{1}{3}-\nu_{1}\right||c_{1}|^{4}+\frac{(1-\alpha)(3-\alpha)}{(2-\alpha)^{2}}|c_{2}|^{2}\bigg).
\end{split}
\ee
For the function $\omega(z)=c_{1}z+c_{2}z^{2}+...$ (with $|\omega(z)|<1, \,z\in\ID$ ) the next relations is valid (see, for example \cite[p.128, expression (13)]{Prokhorov-1984}):
\be\label{eq-24}
|c_{1}|\leq1,\,\, |c_{2}|\leq 1-|c_{1}|^{2},\,\,
|c_{3}(1-|c_{1}|^{2})+\overline{c_{1}}c_{2}^{2}|\leq (1-|c_{1}|^{2})^{2}-|c_{2}|^{2}.
\ee
We may suppose that $a_{2}\geq 0$, which implies that $c_{1}\geq 0$ and instead of relations \eqref{eq-24}
we have the next relations
\be\label{eq-25}
0\leq c_{1}\leq1,\,\, |c_{2}|\leq 1-c_{1}^{2},\,\,
|c_{3}|\leq 1-c_{1}^{2}-\frac{|c_{2}|^{2}}{1+c_{1}}.
\ee
By using \eqref{eq-25} for $c_{1}$ and $c_{3}$, from \eqref{eq-23} we have
\be\label{eq-26}
\begin{split}
|H_{2}(2)| &\leq \frac{4\gamma^{2}}{(1-\alpha)(3-\alpha)}\left[c_{1}(1-c_{1}^{2})+
\left(\frac{(1-\alpha)(3-\alpha)}{(2-\alpha)^{2}}-\frac{c_{1}}{1+c_{1}}\right)|c_{2}|^{2}\right.\\
&\left.+ \mu_{1}c_{1}^{2}|c_{2}| +\left|\frac{1}{3}-\nu_{1}\right| c_{1}^{4}\right].
\end{split}
\ee

Since for $0<\alpha <2-\sqrt{2}$ we have $\frac{(1-\alpha)(3-\alpha)}{(2-\alpha)^{2}}\geq
\frac{1}{2}\geq \frac{c_{1}}{1+c_{1}}$, then by using $|c_{2}|\leq 1-c_{1}^{2}$, from \eqref{eq-26}
after some calculations we obtain

\be\label{eq-27}
|H_{2}(2)|\leq\frac{4\gamma^{2}}{(1-\alpha)(3-\alpha)}F(c_{1}),
\ee
where
\be\label{eq-28}
F(c_{1})=\frac{(1-\alpha)(3-\alpha)}{(2-\alpha)^{2}}+Ac_{1}^{2}+Bc_{1}^{4},
\ee
where
\[
A=\frac{2\gamma-(\alpha^{2}-4\alpha +2)}{(2-\alpha)^{2}}, \\
B=\left|\frac{1}{3}-\nu_{1}\right|-\frac{2\gamma+1}{(2-\alpha)^{2}}.
\]

Further, by using the assumptions of the theorem that
$0<\alpha< 2-\sqrt{2}$ and $0<\gamma \leq\frac{1}{2}(\alpha^{2}-4\alpha +2),$
we easily conclude that $A \leq0 $, while

$$0<\nu_{1}=\frac{(\alpha^{2}-10\alpha +13)\gamma^{2}}{3(1-\alpha)^{2}(2-\alpha)^{2}}
\leq\frac{(\alpha^{2}-10\alpha +13)(\alpha^{2}-4\alpha +2)^{2}}{12(1-\alpha)^{2}(2-\alpha)^{2}}<\frac{13}{12}.$$

If we have that $B\leq 0$, then from \eqref{eq-28} we obtain that
$$F(c_{1})\leq\frac{(1-\alpha)(3-\alpha)}{(2-\alpha)^{2}},$$
and if $B>0$, then
$$F(c_{1})\leq \max\{F(0), F(1)\}=\max \left\{\frac{(1-\alpha)(3-\alpha)}{(2-\alpha)^{2}}
,\left|\frac{1}{3}-\nu_{1}\right|\right\}=\frac{(1-\alpha)(3-\alpha)}{(2-\alpha)^{2}},$$
since
\be\label{eq-30}
\frac{(1-\alpha)(3-\alpha)}{(2-\alpha)^{2}}>\left|\frac{1}{3}-\nu_{1}\right|
\ee
when
$0<\alpha< 2-\sqrt{2}$ and $0<\gamma \leq\frac{1}{2}(\alpha^{2}-4\alpha +2)$ (proven later).
It means that in both  cases we have that
$$F(c_{1})\leq\frac{(1-\alpha)(3-\alpha)}{(2-\alpha)^{2}},$$
which by  \eqref{eq-27} implies
$$  |H_{2}(2)|\leq \left(\frac{2\gamma}{2-\alpha}\right)^{2}.$$

We need to prove the inequality \eqref{eq-30} for appropriate $\alpha$ and $\gamma$.
First, if $\frac{1}{3}-\nu\leq 0$, i.e. if  $0<\nu_{1} \leq\frac{1}{3}$, then , since
$0<\alpha< 2-\sqrt{2}$, we have
$$\frac{(1-\alpha)(3-\alpha)}{(2-\alpha)^{2}}>\frac{1}{2}>\frac{1}{3}-\nu_{1},$$
which implies that \eqref{eq-30} is true.
In case $\nu_1>\frac{1}{3}$, we have that inequality \eqref{eq-30} is equivalent to
$$\frac{(1-\alpha)(3-\alpha)}{(2-\alpha)^{2}}>\frac{(\alpha^{2}-10\alpha +13)\gamma^{2}}{3(1-\alpha)^{2}(2-\alpha)^{2}}-\frac{1}{3}.$$  The last inequality is equivalent
with
\[
\gamma^{2}< \frac{(1-\alpha)^{2}(4\alpha^{2}-16\alpha +13)}{\alpha^{2}-10\alpha +13}.
\]
Since for $0<\alpha< 2-\sqrt{2}$ we have $\gamma \leq\frac{1}{2}(\alpha^{2}-4\alpha +2),$
then for such $\alpha$ we have
$$\gamma^{2}\leq\frac{1}{4}(\alpha^{2}-4\alpha +2)^{2}$$ and from \eqref{eq-30} it is sufficient to prove that \be\label{eq-32}
\frac{1}{4}(\alpha^{2}-4\alpha +2)^{2}\leq \frac{(1-\alpha)^{2}(4\alpha^{2}-16\alpha +13)}{\alpha^{2}-10\alpha +13}
\ee
for $0<\alpha< 2-\sqrt{2}$.
The inequality \eqref{eq-32} is equivalent to
\be\label{eq-33}
( \phi(\alpha):=)\,4(1-\alpha)^{2}(4\alpha^{2}-16\alpha +13)-
(\alpha^{2}-4\alpha +2)^{2}(\alpha^{2}-10\alpha +13)\geq0,
\ee
where $0<\alpha <2-\sqrt{2}$ .
Let's put $\alpha^{2}-4\alpha +2=t$. Then $0<t<2 $ and $\alpha=2-\sqrt{2+t}$ and from \eqref{eq-32}
we have
\[
\phi_{1}(t):=\phi(2-\sqrt{2+t})= \frac14(2+t)\left[30+19t-t^2-(20+6t)\sqrt{2+t}\right].
\]
The function $\phi_{1}$ is continuous function in the interval $[0,2]$. It is easily to check
that
\[ \phi_1'(t) = \frac14\left[ 68+34t-3t^2-(42+15t)\sqrt{2+t} \right]\]
and
\[ \phi_1''(t) = \frac18\left[ 68-12t-45\sqrt{2+t}-\frac{12}{\sqrt{2+t}} \right].\]
iN $\phi_1''$, the second and the third expression reach their minimum on the segment $[0,2]$  for $t=0$, while the last expression for $t=2$. Thus
\[
\phi_1''(t) < \frac18\left( 68-12 \cdot 0 - 45\sqrt{2+0}-\frac{12}{\sqrt{2+2}} \right) = \frac18(62-45\sqrt2) = -0.20\ldots <0,
\]
i.e, $\phi'_1$ is an decreasing function from $\phi'_1(0)=17-10.5\sqrt{2}=2.15\ldots>0$ to
$\phi'_1(2)=-5<0$, which implies that the function $\phi$ attains its maximum in the interval $(0,2)$,
so that
$$ \phi_{1}(t)\geq \min\{\phi_{1}(0),\phi_{1}(2)\}=\min\{15-10\sqrt2, 0\}=0.$$
This means that the inequality given by \eqref{eq-33} is true.

The result of Theorem \ref{18-th-2} is the best possible as the functions $f_{2}$,  defined with
$$\left(\frac{z}{f_{2}(z)}\right)^{1+\alpha}f_{2}'(z)=\left(\frac{1+z^2}{1-z^2}\right)^{\gamma}$$
shows. In this case we have that $c_2=1$, $c_j=0$  when  $j\neq 2$, and consequently, $a_2=a_4=0$, $a_3=\frac{2\gamma}{2-\alpha}$ and $H_2(2)=a_2a_4-a_3^2=-\frac{4\gamma^2}{(2-\alpha)^2}$.
\end{proof}

\label{lastpage}
\end{document}